\documentclass[10pt]{article}
\usepackage{geometry} 
\pdfoutput=1
\usepackage{amsmath,amssymb}  
\usepackage{bm}  

\makeatletter
\newcommand{\Bigint}{\@ifnextchar_\@Bigintsub\@Bigintnosub}
\def\@Bigintsub_#1{\def\@int@subscript{#1}\@ifnextchar^\@Bigintsubsup\@Bigintsubnosup}
\def\@Bigintsubsup^#1{\mathop{\text{\LARGE$\int_{\text{\normalsize$\scriptstyle\@int@subscript$}}^{\text{\normalsize$\scriptstyle#1$}}$}}\nolimits}
\def\@Bigintsubnosup{\mathop{\text{\LARGE$\int_{\text{\normalsize$\scriptstyle\@int@subscript$}}$}}\nolimits}
\def\@Bigintnosub{\@ifnextchar^\@Bigintnosubsup\@Bigintnosubnosup}
\def\@Bigintnosubsup^#1{\mathop{\text{\LARGE$\int^{\text{\normalsize$\scriptstyle#1$}}$}}\nolimits}
\def\@Bigintnosubnosup{\mathop{\text{\LARGE$\int$}}\nolimits}
\makeatother

\makeatletter
\newcommand{\bigint}{\@ifnextchar_\@bigintsub\@bigintnosub}
\def\@bigintsub_#1{\def\@int@subscript{#1}\@ifnextchar^\@bigintsubsup\@bigintsubnosup}
\def\@bigintsubsup^#1{\mathop{\text{\large$\int_{\text{\normalsize$\scriptstyle\@int@subscript$}}^{\text{\normalsize$\scriptstyle#1$}}$}}\nolimits}
\def\@bigintsubnosup{\mathop{\text{\large$\int_{\text{\normalsize$\scriptstyle\@int@subscript$}}$}}\nolimits}
\def\@bigintnosub{\@ifnextchar^\@bigintnosubsup\@bigintnosubnosup}
\def\@bigintnosubsup^#1{\mathop{\text{\large$\int^{\text{\normalsize$\scriptstyle#1$}}$}}\nolimits}
\def\@bigintnosubnosup{\mathop{\text{\large$\int$}}\nolimits}
\makeatother

\begin{document}



\title{
Asymptotic near-efficiency of the ``Gibbs-energy (GE) and empirical-variance" estimating  functions for fitting
 Mat\' ern models -- II: Accounting for measurement errors via ``conditional GE mean''
 \footnote{
 The previous version (version 2) considered both the case with measurement errors (also called ``nugget-effect" or simply ``noise") and the no-noise case. The no-noise case is now in Girard (2016) with more detailed proofs and two additional (w.r.t.~version 2)  results: a consistency result is proved and the restriction  $ \nu \geq 1/2$ is eliminated. This version 3 is devoted to the case with measurement errors, and also gives the analogs of these two additional results.
 }
}


\author{Didier A. Girard}

\date{CNRS and   Grenoble-Alpes University, Lab. LJK, F-38000, Grenoble, France}

\maketitle
\def\vZ  {{\bf z}}

\begin{abstract}
Consider one realization of a  continuous-time Gaussian process $Z$ which belongs to the Mat\' ern family with known ``regularity'' index $\nu >0$. 
For estimating the autocorrelation-range and the variance of $Z$ from  $n$ observations on a 
fine  grid, 
we studied in Girard (2016) the GE-EV method which
simply retains the empirical variance (EV)
 and  equates it to a candidate ``Gibbs energy  (GE)",
i.e.~the quadratic form ${\bf z}^T R^{-1} {\bf z}/n$ where ${\bf z}$ is the vector of  
observations and
$R$ is the autocorrelation matrix for ${\bf z}$ associated with  a candidate range. The present study considers the case where the observation is  ${\bf z}$  plus a Gaussian white noise whose variance is known. We propose to simply bias-correct EV and to replace GE by its conditional mean given the observation.
We show that the ratio of the large-$n$ mean squared error of the resulting 
CGEM-EV
 estimate of the range-parameter to the one of its maximum likelihood estimate, and the analog ratio for the variance-parameter, have the same behavior than in the no-noise case: they both converge, when the grid-step tends to $0$, toward a constant, only function of $\nu$, surprisingly close to $1$ provided $\nu$ is not too large.
We also obtain,  for all $\nu$,  convergence to 1 of the analog ratio for the microergodic-parameter.

\end{abstract}
%
%
%
%


\def\magn   {{b}}
 \def\vZ  {{\bf z}}
 \def\vy  {{\bf y}}
 \def\vw  {{\bf w}}
  \def\vectw  {{\bf w}}
  \def\vecty  {{\bf y}}

 \def\vectz  {{\bf z}}
 \def\Rmat {{R_{\theta}}}
  \def\Rmat0 {{R_{\theta_0}}}

\def\intTorus  {{\int_{-\pi}^\pi}}
\def\varNoise  {{2 \pi}}
\def\varNoise  {{\sigma^2 }}

  \def\matA {{A_{\theta}(\magn)  }} 
  \def\R  {{\rm R}}
 
 \def\defby  {{:=}}

	\def\magn   {{b}}

	\def\vZ  {{\bf z}}
	\def\vz  {{\bf z}}

	\def\vy  {{\bf y}}

	\def\vecty  {{\bf y}}

	\def\matA {{A_{\magn,\theta} }}

\def\matAsquared {{A^2_{\magn,\theta} }} 

	\def\vw    {{\bf w}}

	\def\vu    {{\bf u}}

	\def\R    {{\rm R}}

	\def\vx    {{\bf x}}

	\def\fcontinuous    {{f^*_{\nu,\magn,\theta}Ê}}
		\def\fdelta    {{f^\delta_{\nu,\magn,\theta}Ê}}
	\def\fdeltasigma    {{f^\delta_{\nu,\magn,\theta,\sigmaN}Ê}}

		\def\fdeltaZero    {{f^\delta_{\nu,\magn_0,\theta_0}Ê}}

		\def\gcontinuous {{g^*_{\nu,\theta}Ê}}

\def\gdeltaUnaliased {{g^*_{\nu,\delta\theta}Ê}}
				
	\def\hcontinuous    {{h^*_{\nu,\magn,XXXX\theta}Ê}}

		\def\gdelta    {{g^\delta_{\nu,\theta}Ê}}

\def\hdelta    {{h^\delta_{\nu,\theta}Ê}}

	\def \filter    {{a^\delta_{\magn,\theta}Ê}}
	
	\def \newfilter    {{a_{\delta,\magn,\theta}}}
\def \newNewfilter  {{a_{\magn, \delta \theta}}}

	\def\Rmat0 {{R_{\theta_0}}} 

\def\Rmat {{R_{\theta}}} 

\def\O{{\rm O}}

\def\Cnu{{C_\nu}}
\def\CnuP1{{C_{\nu+1}}}
\def\C0{{C_0}}
\def \const{{C}}

\def\sigmaN  {{\sigma_{\rm N}}}
\def\varNoise  {{\sigma^2_{\rm N}}} 
\def\inVarNoise  {{\sigma^{-2}_{\rm N}}} 

%
 	\section{Introduction}

We consider  time-series of length $n$ obtained by observing, at $n$ equispaced times, a continuous-time process $Z$ which is Gaussian, has mean zero and an autocorrelation function which belongs to the Mat\' ern family with ``regularity'' index $\nu >0$. 
See the Introduction of Girard (2016) and the references therein for comments on this popular family.
%
				We just recall, for notational completeness, that a 
				Mat\'ern processes on  $\mathbb{R}$ can be specified by
				its spectral density over $(-\infty, + \infty )$  where $\theta$ designates the so-called  ``(inverse) range parameter":

		$$ \fcontinuous(\omega)=\tau^2 \,  \gcontinuous(\omega), \,\,\,\, {\rm with}\,\,\,\, \gcontinuous(\omega) := {{ \Cnu   \,  \theta^{2 \nu}}  \over 
				{(\theta^2 + \omega^2)^{\nu+{1\over 2}}}}
				\,\,\,\, {\rm where}\,\,\,\, 
				\Cnu=\frac{\Gamma \left(\nu +\frac{1}{2}\right)}{\sqrt{\pi } \Gamma (\nu )}
				.  \eqno(1.1)$$
%
%
In this paper  $\tau^2$ is the variance of $Z(t)$  (it is easily checked that $\int _{-\infty }^{\infty }\gcontinuous(\omega ) d\omega =1$).
%

As in Girard (2016), we are  concerned here with  ``dense'' grid for the observation times   (or  ``locations'')
 in the sense that the sampling period $\delta>0$  is ``small'' enough. 
Stein (1999, Chapter~3)  shows that a standard (i.e.~fixed $\delta>0$) large-$n$ asymptotic analysis followed by a less standard small-$\delta$ analysis  
yields
useful theoretical insights.
This is precisely the asymptotic framework we  use here.

%
%
%
%

		But, we assume now that 
		there are Gaussian i.i.d.~measurement errors, or, equivalently for the parametric inference point of view we take here,  there is a geostatistical  ``nugget effect", with known variance $\varNoise$.
		And 		we  assume that $\nu$ is known. 
		That is, given  known  $\nu>0,  \delta>0$ and $\sigmaN >0$,  one observes only
a vector  of size $n$ which, after scaling by $\sigmaN$, has a distribution satisfying the model:
		$$
		 \frac{\vecty}{\sigmaN}     \sim  N(0,\magn_0  R_{\theta_0}  + I_n) 
		\,\, {\rm where}\,\, 
 \magn_0 =\frac{\tau_0^2}{\varNoise}       \eqno(1.2)$$
with $I_n$ denoting the identity matrix and
 $R_\theta$ 
		the Toeplitz matrix 
		of coefficients
$ \left[R_\theta\right]_{j,k}= K_{\nu,\theta}(\delta {Ê\vert j-k\vertÊ}), \,\,j,k=1,\cdots,n,$ 
%
with	
$K_{\nu,\theta}(t) =\int _{-\infty }^{\infty }\gcontinuous(\omega ) e^{ i \omega t}d\omega		 $
(see  e.g.~Stein (1999, Section 2.5)  for expressions for these autocorrelation functions $K_ {\nu,\theta}(\cdot)$).
We  can  thus  call  $\magn_0$ the true signal-to-noise ratio (SNR). Notice that one may already expect  that the results of our present study  for the particular case $\magn_0 \gg1$ and, say, $\varNoise=1$,  will approximately coincide with those of the  ``no-noise" situation of Girard (2016) 
(where $\magn_0$  designated the true variance of $Z$).

The CGEM-EV method, introduced in the  first arXiv version of Girard (2012)  and  that we study here,
 is an extension of GE-EV (which was studied in Girard (2016)) to 
 the case of noisy observations  (or nugget-effect) of known variance, that we consider as a ``natural'' extension.
Indeed, recalling that,  firstly GE-EV reverses the roles played by the variance and the range-parameter  in the well known hybrid method of 	Zhang and Zimmerman (2007) (where a ``rough''  estimate
of the range is used) and uses the  ``rough'' empirical variance,
 it seems natural to merely correct 
 this naive, yet near-efficient (in the sense stated by  Girard, 2016), estimator of $\tau_0^2$,  by its known bias.
Thus we  define
 $$\hat\tau_{\rm EV\vert \sigmaN} ^2:= \frac{1}{n}  \vecty^T \vecty - \varNoise
 \,\, \,\, \,\, {\rm and }\ \,\, \,\, \,\, 
 \hat\magn_{\rm EV\vert \sigmaN}  := \frac{\hat\tau_{\rm EV\vert \sigmaN} ^2}{\varNoise}     
 . \eqno{(1.3)}$$ 
 The  second ingredient of CGEM-EV  consists of replacing  the maximization of the likelihood (ML) w.r.t.~$\theta$ by the   following estimating equation in $\theta$:
 denoting by  $A_{\magn,\theta}$  the ``signal extraction'' matrix (see Section 2)
 $$
 A_{\magn,\theta} :=\magn \Rmat {\left( I_n + \magn \Rmat \right)} ^{-1},
 $$
 find, with  $\magn$ fixed at $ \hat\magn_{\rm EV\vert \sigmaN} $    a root $\theta$ of
$$ {\rm CGEM}(\magn,\theta) - \magn \varNoise
\,\,\,\,\,\, {\rm with}\,\, 
 {{\rm CGEM}(\magn,\theta)} :=  \frac{1}{n}   \vecty^T \matA  \Rmat^{-1} \matA \vecty + \left( \frac{\magn}{n}  {\rm tr} {\left( I_n - \matA \right)} \right)  \varNoise.
 \eqno{(1.4)}$$ 
Recall that, if the un-noisy discretely sampled process, say $\vectz$, were observed,  the second ingredient of the GE-EV method (the equation which replaces (1.4))
 would consist of finding the matching between the variance $\tau^2$
 and $n^{-1}  \vectz^T \Rmat^{-1} \vectz$, a quantity we call the  
candidate ``Gibbs energy" (GE, in short) of  $\vectz$.
 On the other hand,
 by classic manipulations, one can check
  that $ {\rm CGEM}(\magn,\theta)$  is the conditional mean of this GE given $\vecty, \sigmaN$ and the candidate $(\magn,\theta)$.
%
Let us now  combine a well known result about the use of likelihood scores in case of incomplete data (e.g.~Heyde, Section 7.4.1), and the remark recalled in Girard (2016)  that 
$n^{-1}  \vectz^T \Rmat^{-1} \vectz  -  \tau^2$ is, up to a strictly positive deterministic factor, the derivative of the log-likelihood of $\vectz$ w.r.t.~$\magn$. 
We thus deduce that  the proposal (1.4) is, in fact, (and still up to a ($>0$) factor) the {\it likelihood score   w.r.t.}~$\magn$ (and {\it not} $\theta \,$ !)  when only $\vecty$ is observed.
%
%
Thus the first heuristic justification in Girard (2016) could be repeated here, except that the analog of the constrained ML $\magn$-estimator function (i.e., in the nototation of Girard (2016), $\theta \rightsquigarrow
 \hat\magn_{\rm ML}(\theta)$; recall that the heuristic given there, was that adjusting $\theta$ so that  this function be matched to a rough estimate $\magn_1$ for the variance is a useful idea, at least in the infill framework) is no more explicit, and the score equation may be unsufficiant to define such a function
 (note also that the theoretical result of Kaufman and Shaby (2013)  only deals with the no-noise case).

In the following,  we  denote by $\hat\theta_ {\rm GEV\vert \sigmaN}$  this  range-parameter estimate
 (in practice, for a reason suggested at the third paragraph of Section 4, we  chose the smallest root in case of multiple roots).

Note that, $\magn$ being fixed at $\hat\magn_{\rm EV\vert \sigmaN}$, computing ${\rm CGEM}( \magn , \theta)$ at candidate $\theta$ does not require to apply $\Rmat^{-1}$ (since, obviously $ \matA  \Rmat^{-1} \matA = \magn  \Rmat \left(I_n +  \magn \Rmat \right)^{-1} $)
and it is thus the condition number of $ \Rmat+  \magn^{-1} I_n $, not of $ \Rmat$ (as it was the case for GE-EV), which controls the numerical stability of the computation. 
Numerical experiments by Lim et al.~(2017) provide a detailed analysis of this condition-number  for the Mat\' ern covariance.
Thus, as it was already known from experiments in kriging or ML computations, numerical instability can be alleviated by adding in the model an, even small, nugget effect. That is, even for un-noisy observations,
 it may be useful (and sometimes mandatory) to use CGEM-EV, instead of GE-EV, with a small a priori fixed $\varNoise$ (the 
impact of such a prior value is studied in the last experiment of Girard (2017)). 
Let us add that in the latter paper, other comments are given (especially in its Secttion~2.3)  about  computational aspects of CGEM-EV as an alternative to ML when $\varNoise$ is known. 
		In  this article, we shall provide
		an asymptotic justification for CGEM-EV, as compared to ML, identical to that already obtained for GE-EV in the no-noise case, 
		except we do not give a precise meaning of the ``small-ness'' of $\delta$ which is sufficient for guaranteeing an asymptotic consistency of $\hat\theta_ {\rm GEV\vert \sigmaN}$.
Recall that the ``$\nu$ not too-far from  $1/2$'' condition  is required to obtain appealing near-efficiency results
 (more precisely, e.g.,~$0< \nu \leq3$ implies a mean-squared-error inefficiency less than $1.33$ in the asymptotic framework we use). In practice,  $\nu>3$ is rarely used, see e.g.~Stein (1999), Gaetan and Guyon (2010)).
We  hope that the theoretical justification obtained here can be extended to more computationally complex settings. 	Indeed, this approach is clearly not limited to observations on a one dimensional lattice, and is potentially not limited to regular grids (a  weighted version, with Riemann-sum type coefficients, of the empirical variance may then be useful).
Successful experiments with CGEM-EV and its Riemann-sum version, with various simulated two-dimensional Mat\' ern random fields, 
are described in Girard~(2017).    
See also the Mathematica Demo  (Girard 2014)   we produced so that 
  any one can
  easily assess  
  CGEM-EV for the case $\nu=1/2$.

	The rest of this article is structured exactly as Girard~(2016),
	except that, in addition, the infill framework  is somewhat discussed at the end of Section~4.

 \section{ Further notations and some properties 
 of the  spectral densities for  Mat\' ern time-series}

Let us recall that, as in Girard (2016),
		we choose the vocabulary here (i.e.~``time"  in place of ``space") since we use in numerous places of the paper the now classical time-series theory.
Set-up (1.2) is 
equivalent to assuming that only a Gaussian time-series  $Z_\delta$, defined by $Z_\delta(i):= Z(\delta i)$, perturbed by a Gaussian white noise, independant of $Z$, is observed at  $i=1,2,\cdots,n$.
From the well known aliasing formula (e.g.~Section 3.6 of Stein (1999)), the spectral density on $(-\pi,\pi]$ of the observed series is 
		                        $$ \fdeltasigma = \varNoise  \left( \magn \,  \gdelta + \frac{1}{2 \pi}\right) \,\,\,\,\,\,{\rm with} \,\,\,\,\
		                         \gdelta (\cdot) := {1 \over  \delta} \sum_{k=-\infty} ^\infty  
		                         \gcontinuous\left( { \cdot  + 2 k \pi \over \delta} \right).
		                             \eqno(2.1)$$
Recall that,
when $\nu -1/2$ is an integer, then $\gdelta$ coincide with particular ARMA spectral densities with  a {\it constrained} vector of parameters.
In order to simplify the  statement of the results here (and their proofs), it is convenient to introduce 
the following weight function $\filter(\cdot)$ over $(-\pi, \pi]$, that we call the candidate filter for given $(\magn,\theta)$
$$ 
\filter(\cdot)    :=  {   \gdelta(\cdot)     \big /   \left({  \gdelta(\cdot)  + (2 \pi \magn)^{-1}}\right).}
\eqno(2.2)$$
Indeed, as is well known from the  signal extraction  literature, $\filter$ is the frequency response  of the  ``optimal (if $\magn,\theta$ were the true parameters)" convolution  of the perturbed series if it were observed over $\mathbb{Z}$;
see e.g.~Section 4.11 of Shumway and Stoffer 2006 for details,  and Girard (2012) also for related well known ``best extracting'' properties of applying the matrix $\matA $.

A function which 
will play an important role in this article (as it was the case in  Girard (2016))
is  the derivative of $\log (\gdelta(\cdot) )$ w.r.t.~$\theta$;
we just recall that it has the following useful expression: $$\hdelta =  { 2 \nu \over \theta }  \left(    1  -  {  g_{\nu+1,\theta}^\delta  \over     \gdelta }  \right),
\quad {\rm where}\quad   
\hdelta :=  \partial  \log (\gdelta )  /\partial \theta
.     \eqno(2.3)$$

For any   $f :[-\pi,\pi ] \rightarrow \mathbb{R}$, s.t.~$\intTorus w(\lambda) f(\lambda) {\rm d}\lambda \neq 0$, where $w(\cdot) >0$ is a
 weight function (we, in fact, only use $w:=[\filter]^2$) 
we define the weighted  coefficient of variation of $f$ by
$$ 
J_w(f) :=
 {  \left\{
 { 1 \over{ \int w}}
  \Bigint  w \left\vert  f -      { 1 \over {\int w}  }    \int  w f     \right\vert^2  \vphantom{{{{{{{{{{{{a^b}^b}^b}^b}^b}^b}^v}^b}^b}^b}^b}^b}
 \right\}        
  \Bigg /   
 \left(        { 1 \over{ \int w }  }    \int w f  \right)^2}
 = 
 { \frac{ 
 { 1 \over {  \int w }}
  \int  w f^2 
  }      
{ \left(        { 1 \over { \int w }  }    \int w f  \right)^2}
 }
 -1  \,\,.
  \eqno(2.4)$$
%

Above and throughout this paper, ``$\int $'' will  denote  integrals over $[-\pi,\pi ]$. Omitting the indexes $\delta$ and $\nu$, we will also use the notation $g_0$
(resp.~$h_0$) for the function $\gdelta$ (resp.~$\hdelta$) when $\theta=\theta_0$.
%
%
%
%
%
%
$B$ (resp.~$\Theta$) will denote any   compact interval not containing $0$ and such that $b_0$ (resp.~$\theta_0$)  is in the interior of $B$ (resp.~$\Theta$) .

We now collect in the following lemmas (whose proof are postponed to an Appendix)   small-$\delta$  equivalences which will be used to prove the results of the following Sections; they might be of interest also for other studies of the Mat\' ern   time-series plus white noise  model:
  \bigskip

\noindent{\bf Lemma 2.1.} {\it For any $\magn >0,\theta >0,\nu >0$ and $k \in\{1,2\}$, we have as $\delta  \downarrow 0$:
{$$
 \int  {  \left( \frac{ \gdelta}{\filter} \right)^2      }  
 \sim
 \int  {  \left(  \gdelta \right)^2      }  
   \sim \frac{c_{1,\nu}}  {\delta\theta},
\,\,\,\,\,\,\,
 \int  {    [ \filter ]^2  \left( \frac{ g_{\nu+1,\theta}^\delta}{ \gdelta} \right)^k      }  
 \sim
   \int  {  \left( \frac{ g_{\nu+1,\theta}^\delta}{ \gdelta} \right)^k  } 
   \sim 2 \pi  \, {c_{2,\nu}}^k  \, \delta\theta$$ 
%
where the constants $c_{1,\nu}, c_{2,\nu}$ are given in Lemma~2.1 of Girard (2016).

}
%
%
  \bigskip

\noindent{\bf Lemma 2.2.} {\it 
For any $\magn >0,\theta >0,\nu >0, k \in\{1,2\}$ and with ${C_{\nu}}$ defined in (1.1), we have as $\delta  \downarrow 0$:
%
$$
{
{{\int { \left[ \filter(\lambda) \right]^k   }  {\rm d}\lambda }
\sim
{
 {2 \delta ^{\frac{2 \nu }{2 \nu +1}}   {(2 \pi  \Cnu  c)^{\frac{1}{2 \nu +1}  }}  {\Gamma \left(k-\frac{1}{2 \nu +1}\right) \Gamma \left(1+\frac{1}{2 \nu +1}\right)}}
 \,\,\,\, \,\,  {\rm where}\,\, \,\, \,\, 
 c=\magn \theta^{2 \nu}.
}}
}
$$
 }}
   \bigskip
   

Now from the fact that ${\delta^{{2 \nu \over 2 \nu +1}} }$ dominates $\delta$ for any $\nu >0$, and the expression (2.3)  of $\hdelta$, the following corollary is easily obtained (by proving the third stated equivalence before the first one):
   \bigskip
   
\noindent{\bf Corollary 2.3.} {\it For any $\magn >0,\theta >0,\nu >0$ and $k \in\{1,2\}$, we have as $\delta  \downarrow 0$, for the weight function $w =  [ \filter ]^2$:
$$  J_{w} \left(  \hdelta \right)
 \sim
 { 1 \over {\int  [ \filter ]^2 }  }
   \int  {  \left(g_{\nu+1,\theta}^\delta  \over     \gdelta\right)^2  },
   \,\,\,\,\,\,\,
 J_{w} \left( \frac{\gdelta}{ [ \filter ]^2}\right)
 \sim
 {\int  [ \filter ]^2 }   \int  {  \left(  \gdelta \right)^2      } 
 %
$$
and
$$
   \int  {    [ \filter ]^2  \left( \hdelta \right)^k      }  
 \sim
   { \left( 2 \nu \over \theta \right)^k  } \int  {    [ \filter ]^2      }.  
$$
}

 \section{Consistency}

Of course, at fixed $\delta$, $\hat\magn_{\rm EV\vert \sigmaN} $  is a consistent estimator of $\magn_0$  (see also (4.1)).
		We first state asymptotic properties of the  (normalized)  
		 estimating equation 
		$	{{\rm CGEM}(\magn,\theta)}  - \magn   {\varNoise} =0$
		and its partial derivatives,
in particular    for $\delta$ ``small'', whose proof only requires  classical techniques and the third equivalence of Corollary~2.3 (see the comments below):

\bigskip

{\noindent \bf Theorem 3.1.} {\it   1)  We have the following   three convergences  in probability, uniform over $B \times \Theta$,   as $n \rightarrow \infty$ :
	\begin{eqnarray*}
		{ \sigma^{-2}_{\rm N}  {{\rm CGEM}(\magn,\theta)  - \magn}  \over{ \left(  n^{-1} {\rm tr}\matA  \right) \magn  }}
		&=&
		\left(  {{ \vecty^T \matA (I-\matA) \vecty \over{ \varNoise  {\rm tr}\matA  }}  - 1 }  \right)
%
\\
		&\rightarrow&     { 1 \over {\intTorus  \filter(\lambda)}   }
		 \bigint_{-\pi}^{\pi}    
		  [   \filter(\lambda) ]^2 
		    {\left(    {\magn_0 {{g^\delta_{\nu, \theta_0}Ê}}(\lambda)  \over \magn {g^\delta_{\nu, \theta}Ê}(\lambda)}  -1 \right)} d\lambda
		    		\,\, \,\,\,\, =: \phi (\delta,\magn,\theta,\magn_0,\theta_0), \,\, {\rm say,}
\end{eqnarray*}
$$		{ \partial  \over {\partial \magn} }
 \left(
		\frac	{{\rm CGEM}(\magn,\theta)}  {\varNoise} - \magn    
				 \right)
\rightarrow  
 %
		{- 1  \over 2 \pi} 
  \left(
  \int_{-\pi}^{\pi}  { {  [  \filter(\lambda)  ]^2   }    }
   d\lambda
-  2 \bigint_{-\pi}^{\pi}  \left(  [  \filter(\lambda)  ]^2 -  [  \filter(\lambda)  ]^3 \right)    {   \left( {  {\magn_0 g_0(\lambda) }  \over  {\magn   \gdelta(\lambda) } } -1 \right) }
  d\lambda
   \right),
$$

$$
				 { \partial  \over {\partial \theta} }
 			 \frac	{{\rm CGEM}(\magn,\theta)}  {\varNoise}    
\rightarrow 
 %
		{- \magn  \over 2 \pi} 
  \left(
  \int_{-\pi}^{\pi}  { {  [  \filter(\lambda)  ]^2   }   \hdelta(\lambda) }
   d\lambda
+ \bigint_{-\pi}^{\pi}  \left(   2 [  \filter(\lambda)  ]^3 - [  \filter(\lambda)  ]^2 \right)  \hdelta(\lambda)   {   \left( {  {\magn_0 g_0(\lambda) }  \over  {\magn   \gdelta(\lambda) } } -1 \right) }
  d\lambda
   \right).
	$$
%
%
%
%
2)  When $\delta \downarrow 0$, we have
		$$\phi(\delta,\magn,\theta,\magn_0,\theta_0) \rightarrow 	
%
{2 \nu }( {2 \nu +1})^{-1}  \,
 \left( {\magn_0 \theta_0^{2 \nu}  \over  \magn \theta^{2 \nu} }  -1
 \right).
 %
%
%
 %
%
		 $$
		 3)  
There
exists a strictly positive function $\bar{\delta}(\nu,\magn_0, \theta_0)$ such that   $0<\delta   \leq \bar{\delta}(\nu,\magn_0, \theta_0)$ implies that the large-$n$ limit in probability of ${ \partial  \over {\partial \theta} }
{\rm CGEM}(\magn,\theta)$
 evaluated at $(\magn_0,\theta_0)$ is stricly negative.
}
\bigskip
			As  it was the case for Part~1 of Theorem~3.1 of Girard (2016), the first part here is in fact not restricted to the Mat\' ern family. Indeed, it only requires regularity conditions on $g^\delta_{\nu, \theta}Ê(\cdot)$, and its strict positivity, which are well fulfilled; and the three limits of 1) can directly be obtained, albeit more tediously than in the no-noise case, from
			 classical large-$n$ theoretical results
			 about quadratic forms constructed from a product of powers, possibly negative, of Toeplitz matrices
			 (e.g.~Azencott and Dacunha-Castelle (1986)).
	
		The second part of Theorem~3.1 is, on the contrary, a consequence of specific properties of the Mat\' ern family, and, in fact, it can be proved by the same techniques as those used
		in Section 3 of  Stein (1999). 
		Let us comment this small--$\delta$ equivalent associated with the first $p$-limit of 1).
		Firstly,
		by examining the analog previous results in the no-noise case, 
		we see that the  first of these previous results
	 is well a ``particular case'' of the first limit above by setting  $\filter$ to $1$,
		 which is well in agreement with  the guess that the no-noise case corresponds to $\filter =1$
		(notice that a similar remark can be made for the terms of the Jacobian given in Proof of Part~1 of Theorem 4.1 of Girard (2016) which are seen as particular 
		case of the second and third limits above with $\filter$ set to $1$).
		Secondly, always compared to the no-noise case,
		   the small-$\delta$ limit of the $p$-limit  (after normalization) of the equation $ \sigma^{-2}_{\rm N}  {{\rm CGEM}(\magn,\theta) = \magn} $ to be solved in $\theta$,
		    is unchanged except for a constant factor (in fact this factor could have been eliminated if we had normalized by $( n^{-1}{\rm tr}\matAsquared) \magn$   in place of  $( n^{-1}{\rm tr}\matA) \magn$,
		  but this is unimportant 
		    and it seems more natural to choose ${\rm tr}\matA$ since it yields a
		    simple expression  for the left-hand term of the first result in Theorem~3.1 (expression given in parentheses)).
		    
		  Let us thus recall that, if $\magn$ is fixed at any value $\magn _1$, then  the unique root $\theta_1$ of this small-$\delta$-large-$n$  equivalent equation  will satisfy $\magn_1 \theta_1^{2 \nu}= \magn_0 \theta_0^{2 \nu}$. 
		This indeed gives some  support to the 
		extension to CGEM-EV of the first heuristic for GE-EV in Girard (2016), as is discussed in the Introduction.
%
%
%

As to the third part of Theorem~3.1, the existence of such a  function $\bar{\delta}(\nu,\magn_0, \theta_0)$ is of course a consequence of the third equivalence of Corollary~2.3 and the strict positivity of $\filter$ (consequence of its definition), since, from the third result of Part 1, 
the limit in probability of ${ \partial  \over {\partial \theta} }
{\rm CGEM}(\magn,\theta)$
 evaluated at $(\magn_0,\theta_0)$ clearly reduces to $- \magn_0 \varNoise (2  \pi)^{-1} \int {  a_0^2 h_0}$ (indeed the second integral vanishes).

   \bigskip
\noindent{\bf Remark 3.2.} 
Since we do not give an explicit form for the upper-bound $\bar{\delta}(\nu,\magn_0, \theta_0)$, this  is not a result as strong as the analog Part 3 in Girard (2016). Anyway, we believe that
even the result in the no-noise case could be improved and we conjecture that the ``local well-posedness'' of the estimating-equation around $\theta_0$ (namely a garantee
 that this derivative at $(\magn_0, \theta_0)$ converges in probability toward a  non-zero value, as $n \rightarrow \infty$) 
 does not requires that  $\delta$ be small.
   \bigskip
   
A Cramer-type consistency can now be proved (as detailed in the Appendix) by using Kessler et al.~(2012) (where
a survey of general asymptotic results for estimating equations is given, see their Section 1.10); precisely:


   \bigskip
   
\noindent{\bf Theorem 3.3.} {\it 
Assume that   $\delta$ is not greater than $ \bar{\delta}(\nu,\magn_0, \theta_0) $, then
there exists a sequence  of roots $  \hat\theta_ {\rm GEV\vert \sigmaN} $ of the CGEM-EV equation (i.e.~(1.4) with $\magn$ fixed at $\hat\magn_{\rm EV\vert \sigmaN}$), as $n$ increases, which converges in probability to $\theta_0$.
} 
\bigskip

\section{\bf  Mean squared error inefficiencies of CGEM-EV to ML for  the variance, range and  microergodic  parameters}

%

As is common in classical (in the sense that  the sampling period $\delta$ is fixed) time-series theory,
 the term ``asymptotic variance of an estimator'',
 denoted ${\rm avar}(\cdot)$, 
 will designate in this paper the variance of the limiting distribution of $\sqrt{n}$ times the error of this estimator; the large-$n$ mean squared error (MSE) of this estimator will refer to $n^{-1}$ times its asymptotic variance.
 As noticed in Girard (2016), we could   consider  a size of  $\lfloor n/\delta \rfloor$   for the $n$-th data set :  this would only  multiply all the asymptotic variances by  $\delta$ and the following near-efficiency statements would be inchanged (see e.g.~Brockwell et al.~(2007)).
\smallskip

    Consider first a simplified setting: the case where the microergodic parameter $c_0=  \magn_0   \theta_ 0^{2 \nu} $ is assumed to be known. (Note that it might be more natural to call ``microergodic parameter'' the product $ \tau_0^2   \theta_ 0^{2 \nu}$   since one may prefer that this parameter does not change with $\sigmaN$; however since   $\sigmaN$ is assumed known, choosing between these two definitions will have no impact on the properties of considered estimators, identical up a known factor).
    
  This assumption of a known $c_0$  is of course restrictive and the following Theorem 4.0 may be thought  of as one of weak practical interest. 
However it is known that ``reasonably accurate'', even if not fully efficient, estimates of $c_0$ can be computed by less expensive approaches than 
ML in numerous contexts, and one could thus condition the model with such a  ``reasonable'' value of $c_0$ plugged-in. 
Recall that one of these possibly reasonable approaches is to  fix $\theta$ at  a prior choice $\theta_1$, 
and to maximize the likelihood only with respect to  $\magn$: in certain common settings, this furnishes reasonable estimates of $c_0$ provided
  the a priori chosen range (i.e.~$\theta_1^{-1}$) is  ``fixed at something larger than the true value ($\theta_0^{-1}$)'', as said in the Section 3.1 of
Kaufman and Shaby (2013)    where an empirical study  well demonstrates this, in the no-noise case; 
and it is expected that this still hold under a noise of known variance, for which case this approach can be straightfuly    extended.

  In this simplified setting, one can equivalently focus either on the estimation of $\magn_0$ by the nonparametric estimate $\hat\magn_{\rm EV\vert \sigmaN}$ defined by (1.4), or that of $\theta_0$  by 
$ \left({c_0   /  \hat\magn_{\rm EV\vert \sigmaN}}\right)^{1 / ({2 \nu})}$. Let us choose the former since the asymptotic limiting law of   $\hat\magn_{\rm EV\vert \sigmaN}$  
has a simple expression (see  e.g.~Azencott and Dacunha-Castelle (1986)), for $\delta$ fixed :
$$
n^{1/2}\left(   \hat\magn_{\rm EV\vert \sigmaN}    -  b_0 \right)
\xrightarrow[]{\cal D} N\left( 0,    4 \pi     {\mathrm{v}_1}  \right)
		\,\,\,\,\,\,\, {\rm as} \,\,\,\,\,  n \rightarrow \infty, 
\,\,\,\,\,\, {\rm where} \,\,\,\,\,
{\mathrm{v}_1} := 
   \magn_0^2  \int {a_0^{-2} g_0^2}
.  \eqno(4.1)$$
Note that the variance ${\mathrm{v}_1}$ used in Girard (2016) does not designate the ${\mathrm{v}_1}$ used above, but this previous $v_1$  is clearly a ``particular case'' by substituting $1$ for $a_0(\cdot)$. 
In fact, the ratio of the present $v_1$ to the previous one,  decreases to $1$   as 
$\magn_0 \rightarrow \infty$,
 because,
 at each $\lambda$, the filter function $a_0(\lambda)$ obviously increases toward its limit $1$, as $\magn_0 \rightarrow \infty$.
 
 Now, by
considering the spectral density model $f(\magn, \theta) : \lambda \rightsquigarrow \sigmaN^2 \left( \magn  {g^\delta_{\nu, \theta}Ê}(\lambda) + (2 \pi)^{-1}  \right)$ 
with $\magn \theta^{2 \nu}=c_0$, 
as a function  of  only $\magn$, easily 
establishing that ${\partial \log ( f(c_0/\theta^{2 \nu} , \theta)) \over  \partial \theta } =  \filter (\hdelta - {{2 \nu}  \theta^{-1}} )$  
where $\magn =c_0/\theta^{2 \nu}$
and 
using that ${\partial (c_0 / \magn)^{1/ (2 \nu)} \over  \partial \magn}  = {(2 \nu)}^{-1} (\theta / \magn)$ where $\theta = (c_0 / \magn)^{1/ (2 \nu)}$,   the asymptotic Fisher information w.r.t.~$\magn$ is deduced 
and is seen to be $>0$ (from  
the expression (2.3) of $h_0$ and the fact $\filter>0$).
%
Thus, by an application  (similar, but  easier here,  to the way of establishing (4.3) below) of now classical time-series theory (e.g.~Azencott and Dacunha-Castelle (1986))
one  obtains for the ML 
maximizer  over 
$B$ under  $\magn \theta^{2 \nu}=c_0$, now denoted $\hat
%
\magn_{\rm ML\vert c_0,\sigmaN}$,   as  $n \rightarrow \infty$:
  $$
n^{1/2}\left(   \hat\magn_{\rm ML\vert c_0,\sigmaN}    -  b_0 \right)
\xrightarrow[]{\cal D}  N\left( 0,         {\rm avar}(\hat\magn_{\rm ML\vert c_0,\sigmaN})   \right),
%
$$ where $$
{\rm avar}(\hat\magn_{\rm ML\vert c_0,\sigmaN}) 
:= 4 \pi
b_0^2    \left( {{{2 \nu}  \theta_0^{-1}}}\right)^2 
  %
  {\left( \int {   a_0^{2}     \left( h_0- {{2 \nu}  \theta_0^{-1}}\right)^2  }\right)^{-1} } .  
 \eqno(4.2)$$
%
%

   \bigskip
   
Now by using,  in  (4.2), the expression (2.3) of $h_0$, and the first and second
   equivalences of Lemma 2.1, one obtains:

\noindent{\bf Theorem 4.0.} {\it 
The large-$n$ MSE 
  inefficiency of $\hat\magn_{\rm EV\vert \sigmaN}$  relative to the ML estimator of the SNR $b_0$,
   when $c_0=  \magn_0   \theta_ 0^{2 \nu} $ is known, i.e.~$I^0_{\delta,\magn_0,\theta_0}Ê := {  {4 \pi \mathrm{v}_1} / {\rm avar}(\hat\magn_{\rm ML\vert c_0,\sigmaN}) }$, 
  satisfies  (with $C_{\nu}$  defined in (1.1)):
  $$ 
  I^0_{\delta,\magn_0,\theta_0} \rightarrow 
     \frac{{C_{\nu+1}}^2}{{C_{2\nu+1/2}}  {C_{3/2}}}   =
{
 \frac{\sqrt{\pi }}{2 }   \left(\frac{\Gamma \left(\nu + {3/2}\right)}{\Gamma \left(\nu + {1}\right)}\right)^2 
 }
 \Bigg /
 {
  \frac{\Gamma \left(2\nu + {1}\right)}{\Gamma \left(2\nu + {1/2}\right)}
  }
 =: 
 {\rm ineff}(\nu)
   \,\,\,\,\,\, {\it as} \,\,\,\, \delta \downarrow 0.$$
} 

   It is important to notice that this definition of the constant  $ {\rm ineff}(\nu)$ coincides with the one used in  Girard (2016) for  the no-noise case.
   In the particular case $\nu =1/2$, 
then $ {\rm ineff}(1/2)=1$.
(Note there was a typographical error in Girard (2016) in the second expression of $ {\rm ineff}(\nu)$: precisely the big  fraction slash was omitted;
    however Table~4.1 of  Girard (2016) which displayed
   numerical values of $ {\rm ineff}(\nu)$ for  certain values of $\nu$ is exact.)
The Table~4.1 of  Girard (2016) is not repeated here.   
%
%
%
We only wish to emphasize  that the departure from $1$ of $ {\rm ineff}(\nu)$  as $\nu$ increases, is rather slow.

A second good 
news it that
this inefficiency is not function of the true range $\theta_0$ or the SNR $\magn_0$.
Since it could be expected that these  small-$\delta$-large-$n$-inefficiencies become close to those obtained in the no-noise case only  under  $\magn_0 \gg1$, the result that 
$\magn_0$ has no impact on $ {\rm ineff}(\nu)$ may be thought as rather surprising.  
Recall that the asymptotic inefficiencies in Girard (2016) show the absence of any impact of $\theta_0$ in the no-noise case; this was less surprising because this was already known in the case  $\nu=1/2$ from  the efficiency result of  Kessler (1997) concerning the naive empirical variance.
Thus Theorem~4.0 is a neat extension of the efficiency result of  Kessler (1997)   to the case of measurement noise (after natural  
bias-correction of this empirical variance by subtracting  $\varNoise$),
  and a  ``near-extension" when $\nu$ does not depart too much  from $1/2$ in the sense that 
  $ {\rm ineff}(\nu)$ stays close to $1$, whatever $\magn_0$ may be.

\bigskip

   \smallskip

   %
Now let us return to the case   $\magn_0$ and $\theta_0$  unknown.
   Let $( \hat\magn_{\rm ML\vert \sigmaN}, \hat\theta_ {\rm ML\vert \sigmaN} )$ be a maximizer of the likelihood function over 
$B \times \Theta$ when $\sigmaN$ is known.
One can use arguments similar to those used in Girard (2016) where the asymptotic behavior of the ML estimator was described in the no-noise case:
 the derivation of the  asymptotic information matrix (see Theorem 4.3 of Chapter XIII of Azencott and Dacunha-Castelle (1986))  is classic, albeit more tedious; and the final expressions are relatively simple modifications, by merely adding in appropriate places the weight function  $a_0^2$, precisely: 
  $( \hat\magn_{\rm ML\vert \sigmaN}, \hat\theta_ {\rm ML\vert \sigmaN} )$ is a.s.~consistent and satisfies, as $n \rightarrow \infty$:
    $$
n^{1/2}\left(
\left[
\begin{array}{ccc}
 \hat\magn_{\rm ML\vert \sigmaN}  \\
 \hat\theta_ {\rm ML\vert \sigmaN}
\end{array}
\right]
-
\left[
\begin{array}{ccc}
 \magn_0  \\
 \theta_ 0
\end{array}
\right]
\right)
\xrightarrow[]{\cal D}
N\left( 
\left[
\begin{array}{ccc}
0 \\
0
\end{array}
\right],
4  \pi 
\left[
\begin{array}{cc}
  \sigma_1^2    &   \sigma_{12}\\
 \sigma_{12} &    \sigma_2^2
\end{array}
\right]
\right)
,
$$with  $$
\left[
\begin{array}{cc}
 \sigma_1^2  \\
 \sigma_{12} \\
 \sigma_2^2 \end{array}
 \right]
:=
{ {\left\vert\int { a_0^{2}  h_0}\right\vert^{-2} }  
J_{ a_0^{2}}
(h_0)^{-1}} \left[
\begin{array}{cc}
    \magn_0^2  \int {a_0^{2}  h_0^2}\\
    -  \magn_0     {   {\int {a_0^{2}   h_0}}     }\\
 { \int {  a_0^2}   }
 \end{array}
 \right].  
\eqno(4.3)
$$
  


Again 
note that the
components of the vector $( \sigma_1^2,
 \sigma_{12},
 \sigma_2^2)$ used in Girard (2016) 
 and those of 
 $( \sigma_1^2,
 \sigma_{12},
 \sigma_2^2)$ used above, are
 asymptotic equivalents as 
 $\magn_0 \rightarrow \infty$,
 because  in addition to $a_0 \approx 1$, 
the functional  $J_{a_0^2}(\cdot)$ also becomes close, for large SNR, to the simpler functional $J(\cdot)$ used in Girard (2016).

\smallskip
\smallskip
\smallskip
\smallskip

Concerning CGEM-EV, we claim:

\smallskip
\smallskip
   \smallskip
\noindent{\bf Theorem 4.1.} {\it Assuming $\delta < \bar{\delta}(\nu,\magn_0, \theta_0)$,
let $  \hat\theta_ {\rm GEV\vert \sigmaN} $ be a consistent root of the CGEM-EV equation  (i.e.~(1.4) with $\magn$ fixed at $\hat\magn_{\rm EV\vert \sigmaN}$ and $\varNoise$ is the true noise variance).

\noindent 1)   As $n \rightarrow \infty$
 $$n^{1/2}\left(
\hat\theta_ {\rm GEV\vert \sigmaN}
-
\theta_0
\right)
\xrightarrow[]{\cal D} N\left( 0,   4 \pi  {\mathrm v}_2  \right)
%
\,\,\,\,\,\, {\it where} \,\,
{\mathrm v}_2 =     {  {\left\vert{\int {  a_0^2 h_0}}\right\vert^{-2} } }     J_{a_0^2}(g_0 / a_0^2)     {\int {  a_0^2}} \,\,  .
$$
   %
2) 
The large-$n$ MSE
  inefficiency of CGEM-EV to ML for $\magn_0$ (resp.~for $\theta_0$) 
  being defined by 
$I^1_{\delta,\magn_0,\theta_0}Ê := {  {\mathrm{v}_1} / \sigma_1^2}$ (resp.~$I^2_{\delta,\magn_0,\theta_0}Ê := { {\mathrm{v}_2}   / \sigma_2^2}
=J_{a_0^2}(g_0 / a_0^2) J_{ a_0^{2}}(h_0)$),
%
these two inefficiencies have the following common  limit  (with ${{\rm ineff}(\cdot)}$ as in Theorem~4.0):
%
  $$ 
  I^i_{\delta,\magn_0,\theta_0}Ê  \rightarrow   
 {{\rm ineff}(\nu)}
    \,\,\,\,\,\, {\it as} \,\,\,\, \delta \downarrow 0,\,\,\,
    \,\,\, {\it for} \,\,\,\,\,i\in \{1, 2\}    .
  $$
}
\noindent{\bf Proof:} Part~1 is proved in the Appendix. The limit of both $I^1_{\delta,\magn_0,\theta_0}$ and $I^2_{\delta,\magn_0,\theta_0}$  Êis directly deduced from the equivalences stated in Lemma~2.1 and  Corollary~2.3.

Again, 
as noticed for $\sigma_2^2$ above,
the ${\mathrm{v}_2}$ used in Girard (2016) 
is the limit value of the ${\mathrm{v}_2}$ defined in Part~1 above as $\magn_0 \rightarrow \infty$ for fixed $\delta$.

Thus, as in the case $c_0$ known,  the CGEM-EV estimates of $\magn_0$ and $\theta_0$ are asymptotically nearly efficient provided $\nu$ is not too large, asymptotic full-efficiency being reached for $\nu$ close to 1/2.
Notice it is rather surprising  that these 
 small-$\delta$ large-$n$ inefficiencies are  not function of the underlying
  $\theta_0$ or of the underlying $\magn_0$. 

The remark claimed in  Girard (2016) that the knowledge of $c_0$ does not improve (in terms of small-$\delta$-large-$n$ MSE) 
the  performance of  ML estimation of $\theta_0$ 
or the performance of the alternative to ML we have introduced,
can be also claimed in the present setting of known error variance  (this extension is still also easily checked).

\bigskip
Let us now consider
the estimation of the microergodic parameter 
$c_0$.
By the classical delta-method, one  directly infer from (4.3)  that the asymptotic variance
of
$  \hat c_{\rm ML\vert \sigmaN}
 := \hat\magn_{\rm ML\vert \sigmaN}   \hat\theta_ {\rm ML\vert \sigmaN}^{2 \nu} $
 is
  $4 \pi c_0^2 { {\left\vert \int {  {a_0^2} h_0}\right\vert^{-2} }  J_{a_0^2}(h_0)^{-1}}   
 {\int   \left\vert{  a_0   {  \left( h_0- {{2 \nu} \theta_0^{-1}}\right)}}   \right\vert}^2$ (note that there was a typographical error in Girard (2016) : ``$\Big\vert \int  $'' must be replaced by ``$\int   \Big\vert$''   for the related variance with $a_0\equiv 1$ there).
On the other hand, 
  a similar derivation (albeit more tedious than in the no-noise case) can be done for $ \hat c_ {\rm GEV\vert \sigmaN}$ starting 
  from the asymptotic covariance matrix (detailed in the Appendix) of the vector $(\hat\magn_{\rm EV\vert \sigmaN},  \hat\theta_ {\rm GEV\vert \sigmaN})$ 
  and this gives $ {\mathrm v}_3$ below.
  Now one can easily deduce
  (still using the expression (2.3) of $\hdelta$, Lemma~2.1 and Corollary~2.3)
   that, this time, full-efficiency holds for {\it any}  $\nu>0$, more precisely:
 
\bigskip

\noindent{\bf Theorem~4.2. }{\it Assuming  $\delta < \bar{\delta}(\nu,\magn_0, \theta_0)$,
let $ \hat c_ {\rm GEV\vert \sigmaN} :=  \hat\magn_{\rm EV\vert \sigmaN}  \hat\theta_ {\rm GEV\vert \sigmaN}^{2 \nu}  $ where $  \hat\magn_{\rm EV\vert \sigmaN}$ and $\hat\theta_ {\rm GEV\vert \sigmaN} $ are defined as in Theorem~4.1, we have, with $c_0=  \magn_0   \theta_ 0^{2 \nu} $ :

\noindent{1)  as $n \rightarrow \infty$}
    $$
n^{1/2}\left(    \hat c_ {\rm GEV\vert \sigmaN}  -  c_0 \right)
\xrightarrow[]{\cal D} N\left( 0,   4 \pi  {\mathrm v}_3  \right)
$$ where $$
{\mathrm v}_3 =
{c_0^2 \over {
\smallint a_0^2  \, \,
{\left\vert \smallint { a_0^2  h_0}\right\vert^{2} }
}}
 \left(  J_{\magn_0,\theta_0}(g_0  )   {\left\vert  {\int     a_0^2 \left( h_0- {{2 \nu} \theta_0^{-1}}\right)}   \right\vert}
 ^2
 +
 {\left\vert  {\int  a_0^2 h_0}   \right\vert^2} \right).
%
  $$
2) 
The large-$n$ MSE
  inefficiency of CGEM-EV to ML 
  for $c_0$ is   
$I^3_{\delta,\magn_0,\theta_0}Ê := {4 \pi  {{\mathrm v}_3} / {{\rm avar}(\hat c_{\rm ML\vert \sigmaN} )}}$ and it holds that
  $ 
  I^3_{\delta,\magn_0,\theta_0} \rightarrow 1$ as $\delta \downarrow 0 ;  
 $
 more precisely,
      $$
4 \pi {{\mathrm v}_3}     \sim     {{\rm avar}(\hat c_{\rm ML\vert \sigmaN} )}       
\sim   { 2 c_0^2   \over  (2 \pi)^{-1}  \int {  a_0^2} }   
 \sim   {  2 \pi c_0^2   {{(2 \pi  \Cnu   c_0)^{\frac{-1}{2 \nu +1}  }}} \over  
 {
  {      \Gamma \left(2-\frac{1}{2 \nu +1}\right) \Gamma \left(1+\frac{1}{2 \nu +1}\right)}
  }
  }
   \delta^{\frac{-2 \nu }{2 \nu +1}} 
  {\it    \, \,\,  as  \, \, \, \,} \delta \downarrow 0.
\eqno(4.4)  $$
}


Again, the variance ${\mathrm{v}_3}$ used in Girard (2016)    is clearly obtained from the ${\mathrm{v}_3}$ here by substituting $1$ for $a_0(\cdot)$. 
%
As we remarked for the no-noise case,  this full-efficiency of CGEM-EV concerning $c_0$, even for large $\nu$,  may again be   thought of as a less surprising result  than Part~2 of Theorem~4.1.
Indeed this full-efficiency is  suggested 
by the infill-asymptotic efficiency,  mentionned as a heuristical partial justification in the Introduction,  of 
 our ``proposed" estimator of $c_0/\magn_1$, with $\magn$ fixed at ``any''  $\magn_1$
(we use here quotation marks, only for reminding the sobering fact 
that the strategy of using an arbitrarily fixed $\magn_1$ may provide poor estimates in practice, and  thus we do not actually propose it). 

\bigskip
\noindent{\bf Remark 4.3. }{One can make incidental remarks for the ``case''  $\delta =1/n$.
For  the particular case $\nu=1/2$,
notice that $n^{-1}$ times the right-hand expression of this small-$\delta$ equivalent  (4.4) is  identical, by setting $\delta :=1/n$, to 
$ {4\sqrt{2} {c_0}^{3/2}}{n}^{-1/2}$, 
that is, well coincides with
 the  variance, established in the infill asymptotic framework by Chen et al.~(2000), of the normal approximation of the law of the ML estimator of $c_0$.
 See also Zhang and Zimmerman~(2005) for a detailed rigorous  study of this ``reconciliation'' between the two asymptotic frameworks.
Naturally  one can thus conjecture that a such   coincidence still holds beyond the case $\nu=1/2$, namely, that $n^{-1}$ times 
 the  right-hand expression of  (4.4), with $1/n$ substituted for $\delta$, furnishes the infill asymptotic variance for both the ML  or the CGEM-EV estimator of $c_0$ for any $\nu >0$.
For instance, this would furnish a variance of $({16}/{3})c_0^{7/4} {n}^{-1/4}$  for $\nu=3/2$.
Notice that this latter variance also coincides with a related variance  for the integrated Brownian motion plus  white noise model,
 which could be deduced form the Fisher information given by Theorem~2.3 of Kou (2004) (take $r=2,s=0$ using his notation
 and use that 
${\Gamma \left(2- {(2 \nu +1)}^{-1}\right) \Gamma \left(1+{(2 \nu +1)}^{-1}\right)}   =
{(2 \nu +1)\,  
{\rm B} \left(r- {(2 \nu +1)}^{-1},s+ {(2 \nu +1)}^{-1} \right)} $ where ${\rm B}(\cdot,\cdot)$ is the Beta function, see e.g.~Weisstein (2019)).
}

\bigskip
\noindent{\bf Remark 4.4. }{
Since the framework of 
Zhang and Zimmerman~(2005)
actually does not imposes that  $\sigmaN$ be known,
 their results compared with (4.4) (where $\sigmaN$ is known) show  that by estimating  $\sigmaN$
 (which is often more easy to estimate than $c_0$; see Chen et al.~(2000) for the meaning of such a claim in the infill framework), by the ML principe, 
 one  adds no further error, at least in the  small-$\delta$-large-$n$ framework, to the ML estimator of $c_0$.
 Is is naturally expected that this still holds beyond the case $\nu=1/2$.
 An extension of CGEM-EV to the case $\sigmaN$ unknown, is commented in the following Discussion.}

 \section{\bf  Discussion}


CGEM-EV is thus a natural  extension of GE-EV to the case with measurement errors of known variance,   via bias-correction of the naive empirical variance and  replacement of the unobserved GE function by the  conditional GE mean function. 
We have proved here that identical near-efficiency results still hold, not only in the particular case $\magn_0 \gg1$ for which such a similarity could be expected. 
One may be surprised by the fact that these efficiency results hold for any fixed SNR $\magn_0$, even small.
However, one must keep in mind  that these results deal only with  large-$n$ asymptotics at $\delta$ fixed, always followed  by a small-$\delta$ analysis: one may guess that for  $\magn_0$ too small, very large $n$ may be required to ``see" the stated asymptotic behaviors, and even  increasingly large as  $\delta$ decreases. 
Recall also that, even in the no-noise case, when $\delta$ decreases to $0$, larger data sizes are  required to be able to accurately approximate the actual law of any one of these estimates of $\theta_0$ (or $\magn_0$) by its asymptotic form (indeed this is well known for ML estimates in the case $\nu=1/2$  thoroughly studied  by Zhang and Zimmerman 2005)).
An asymptotic comparison  of CGEM-EV to ML deserves thus a futur study also in the finite-$\delta$ case.


Concerning the problem of estimating multidimensional stationary Mat\'ern fields observed   on a lattice under i.i.d.~noise, as already noticed  for  GE-EV in Girard (2016) in the no-noise case, the  
CGEM-EV approach  is directly applicable, in theory, provided $ \theta$ remains a scalar parameter, e.g.~for isotropic autocorrelations.
We refer to  Girard (2017) for a rather extensive  empirical comparison of CGEM-EV to ML in the two-dimensional case,
with randomized-traces used instead of the exact traces of (1.4).
There, it is noticed, in particular, that difficulties appear in case of  ``too small'' SNRs and ``too  smooth'' fields (for example $\nu  \geq 3/2$) both for CGEM-EV and ML. 
The application of CGEM-EV can be  practical for very large lattice sizes, even with missing data, as soon as computing analogues of $ \left(I_n +  \hat\magn_{\rm EV\vert \sigmaN}  \Rmat \right)^{-1}\vecty$ 
for candidate $\Rmat$, can be done by fast  iterative algorithms; indeed each iteration can be fast since applying $\Rmat$ to a vector can  reduce to three multidimensional discrete (inverse) Fourier transforms.
Thus extensions of the asymptotic results of this paper to such  multidimensional  fitting problems clearly  deserves a detailed study.
And a comparison  of CGEM-EV with the classical (tapered) Whittle-likelihood maximization for such  multidimensional  fitting problems should be of interest.

CGEM-EV might be extended to the important case of unknown noise variance, via, at least,   two simple ways which are described and experimented in Girard (2018). 
One of these two ways is to 
add, as a second simple estimating equation, the first derivative of the likelihood w.r.t. $\varNoise$ for  given $\magn$ and $\theta$ equated to zero. 
Simulations in Girard (2018) demonstrate rather good performance of these two approaches, which thus  warrants further exploration. 

\bigskip
\bigskip
\noindent{\bf REFERENCES } 

{\small

\bigskip \noindent
Azencott, R., Dacunha-Castelle, D., 1986. 
Series of Irregular Observations: Forecasting and Model Building. Springer. New York. 

\bigskip \noindent
Brockwell, P. J., Davis, R. A., and Yang Y., 2007. Continuous-time Gaussian autoregression.
Statistica. Sinica, 17, 63-80

\bigskip \noindent
Chen H.S., Simpson D.G. and Ying Z. (2000). Infill asymptotics for a stochastic process model with measurement error. Statistica Sinica 10, 141-156.
%

\bigskip \noindent
Gaetan, C., Guyon, X., 2010. Spatial Statistics and Modeling. Series: Springer Series in Statistics.  

\bigskip \noindent
Girard, D.A., 2012. 
Asymptotic near-efficiency of the ``Gibbs-energy and empirical-variance'' estimating functions for fitting Mat\' ern models to a dense (noisy) series.
 Preprint version 2 http://arxiv.org/abs/0909.1046

\bigskip \noindent 
Girard, D.A., 2014. "Estimating a Centered Ornstein-Uhlenbeck Process under Measurement Errors".
From  Wolfram Demonstrations Project.

\noindent
http://demonstrations.wolfram.com/EstimatingACenteredOrnsteinUhlenbeckProcessUnderMeasurementE/

\bigskip \noindent
Girard, D.A., 2016. 
Asymptotic near-efficiency of the ``Gibbs-energy and empirical-variance'' estimating functions for fitting Mat\' ern models - I: Densely sampled processes.
Statistics \& Probability Letters, 110, 191-197

\bigskip \noindent 
Girard, D.A., 2017. 
Efficiently estimating some common geostatistical models by `energy-variance matching' or its randomized `conditional-mean' versions. Spatial Statistics, 21(Part A),  1-26.   

\bigskip \noindent
Girard, D.A., 2018. 
"Estimators of a Noisy Centered Ornstein-Uhlenbeck Process and Its Noise Variance".
From  Wolfram Demonstrations Project.

\noindent
http://demonstrations.wolfram.com/EstimatorsOfANoisyCenteredOrnsteinUhlenbeckProcessAndItsNois/

\bigskip \noindent
Heyde, C.C., 1997.
Quasi-Likelihood And Its Application: A General Approach to Optimal Parameter Estimation.   Springer Series in Statistics

\bigskip \noindent
Kaufman, C,  Shaby, B., 2013.  
The role of the range parameter for estimation and prediction in geostatistics.
 Biometrika, 100, 473-484

\bigskip \noindent
Kessler M., 1997. Simple and explicit estimating functions for a discretely observed diffusion
process. Scand. J. Statistics 27, 65-82.

\bigskip \noindent
Kessler, M., Lindner A.,   Sorensen, M., 2012. Statistical Methods for Stochastic Differential Equations.
Section 1.10 General asymptotic results for estimating functions.  CRC Press.

\bigskip \noindent
S.C. Kou, On the efficiency of selection criteria in spline regression, Probab. Theory Relat. Fields 127 (2003), pp. 153-176.

\bigskip \noindent
Lim, C.Y., Chen C.-H., Wu W.-Y., 2017. Numerical instability of calculating inverse of spatial covariance matrices.
    Statistics \& Probability Letters, 
    129, 182-188
    
\bigskip \noindent
Stein, M.L., 1999. Interpolation of Spatial Data: Some Theory for Kriging. Springer.

\bigskip \noindent
Weisstein, E. W., 2019.
"Beta Function." From MathWorld--A Wolfram Web Resource. 

\noindent
http://mathworld.wolfram.com/BetaFunction.html,

\bigskip \noindent
Zhang, H., 2004. Inconsistent estimation and asymptotically equivalent interpolation in model-based geostatistics.  J.~Amer.~Statist.~Assoc.  99, 250-261. 

\bigskip \noindent
Zhang H.  and Zimmerman  D.L., 2005. Toward reconciling two asymptotic frameworks in spatial statistics. {  Biometrika},   { 92}, 921-936. 

\bigskip \noindent
Zhang H., Zimmerman  D.L., 2007. Hybrid estimation of semivariogram parameters. Mathematical Geology,  39(2), 247-260. 


}

\bigskip
\bigskip

\appendix{
\noindent
  	{\large \bf APPENDIX }
}

\bigskip
Let us prove Lemma~2.2 before    Lemma~2.1. This can be done by substituting for $\filter$ an ``un-aliased'' version defined below.
So we first establish the following Lemma which states the order of the differences between exact and un-aliased versions:
\bigskip

\noindent{\bf Lemma A.1.} {\it   For any $\lambda \in [-\pi, \pi]$, any $(\magn,\theta) \in B \times \Theta$, assuming $\delta$ bounded, e.g.~$<1$, 
there exists constants $\underbar c_1,\underbar c_2 >0$,
 $\bar{c}_1,\bar{c}_2 < \infty$  only functions of $\nu$, such that:

\noindent  1)
$$
\underbar c_1  \,  \delta^{2 \nu}  \leq   g_{\nu,\delta\theta}^*(\lambda) \leq   g^\delta_{\nu,\theta} (\lambda)  \leq g_{\nu,\delta\theta}^*(\lambda) +  \bar{c}_1 \delta^{2 \nu},$$ 
 
\noindent 2) 
$$ \underbar c_2  \,  \delta^{2 \nu}  \leq
a_{\magn,\delta\theta}^*(\lambda) \leq   a^\delta_{\magn,\theta} (\lambda)  \leq a_{\magn,\delta\theta}^*(\lambda) + \bar{c}_2 \delta^{2 \nu},
 \,\,\,\, {\rm where}\,\,\,\,
a^*_{\magn,\cdot }   :=  {  \magn g^*_{\nu,\cdot }    \over  { \magn g^*_{\nu,\cdot } + (2 \pi)^{-1}}}.$$

}
\bigskip
\noindent {\bf Proof.}
Part~2 of  Lemma A.1 is a direct consequence of its Part~1 (the second and the third inequalities of Part~2 are
 immediate consequences after noticing that, for the function $w: x \rightsquigarrow  {1}/ \left({1+x^{-1}}\right)$, 
 we have that $0 < x_1 < x_2$  implies 
$0<w(x_2)  - w(x_1)  < x_2- x_1 $;  in the same way, the inequality $w(2 \pi b \underbar c_1  \,  \delta^{2 \nu} )  \leq a_{\magn,\delta\theta}^*(\lambda)$ is a consequence of the first inequality in Part 1; 
it suffice then to use  
$\delta<1$ in $ w(2 \pi b \underbar c_1  \,  \delta^{2 \nu} ) =   {2 \pi b \underbar c_1  \,  \delta^{2 \nu}   }/   \left({ 2 \pi b \underbar c_1  \,  \delta^{2 \nu}   +1}\right)$ 
to obtain the lower bound $\underbar c_2 \delta^{2 \nu} $. So let us prove Part 1. The first inequality and the second one are easily checked.
Let us prove the third one.
Letting $\alpha = \delta\theta$, by arguments identitcal to those used in a step of the proof of Lemma~2.1 in Girard (2016), we can obtain that $\lambda > -\pi$ implies that,
 for any $k\geq 1$, 
 $g_{\nu,\alpha}^*(\lambda+ 2 \pi k)  \leq   g_{\nu,\alpha}^*( 2 \pi (k-1/2))    \leq (\Cnu / (2 \pi)^{2 \nu +1} ) \,  \alpha^{2 \nu} / (k-1/2)^{2 \nu +1}$; and summing these terms over $k=1,2,\cdots$, thus gives a term   $\O (\delta^{2 \nu})$ (notice that $\O (\delta^{2 \nu +1})$ was obtained in the mentioned step since we  considered there $g_{\nu,\theta}^*(\cdot / \delta)$ which can be easily checked to be $\delta g_{\nu,\alpha}^*(\cdot )$). This is shown similarly (except we use  $\lambda < \pi$) for the sum over $k=-1,-2,\cdots$. Combining these two results  gives the claimed bound for
 $\gdelta (\lambda) -  g_{\nu,\alpha}^*(\lambda) =   \sum_{k\neq 0}g_{\nu,\alpha}^*(\lambda+ 2 \pi k)$.

\bigskip

%
%

\bigskip

\noindent{\bf Proof of Lemma~2.1}    Note that the concise expression $c_{2,\nu}^k$ results from $C_{3/2} = 2 \pi {C_{1/2}}^2$.
The first equivalence of  Lemma 2.1 is  trivial by develloping  
$\int  { \left[ { \gdelta}/{\filter} \right]^2 } = \magn^{-2} \int  {\left[ \magn \gdelta + 1/(2 \pi)\right]^2 }    = 
 \int  {  \left[ \gdelta\right]^2 }   + 1/( \pi \magn)    \int  { \gdelta }  +  1/(2 \pi \magn)^2  $  and recalling that $\int  { \gdelta } =1$.
 As to the second equivalence, let us examine the Proof of Lemma~2.1 of Girard (2016).
 We make here  the same change of variable
 $\omega = \lambda / \delta$ (note a typographical error: the fraction slash was omitted in Girard (2016)). The integrand is now modified only by the factor $\filter(\delta \omega)$. Since this factor is $<1$, the fact that 
 Lebesgue dominated convergence theorem is applicable has not to be re-proved, and: 
\begin{eqnarray*} 
\delta^{-1}   \int_{-\pi }^{\pi }   {  \left( \filter(\lambda)  {{g^\delta_{{\nu +1},\theta}(\lambda) }\over {g^\delta_{\nu,\theta}(\lambda) }} \right)^2} d\lambda    
 &=& \int _{-\pi/\delta }^{\pi/\delta }    
\left[  \filter(\delta \omega)   \right]^2		                    { \left( \sum_{k=-\infty} ^\infty    { g^\ast_{{ \nu +1},\theta}\left( { \omega  + 2 k \pi / \delta} \right)} \right)^2 \over
		                        {\left( \sum_{k=-\infty} ^\infty  \gcontinuous\left( { \omega  + 2 k \pi / \delta} \right) \right)^2}  }  d\omega \\
&\rightarrow&  \int _{-\infty }^{\infty }    
		                   \left(  {    { g^\ast_{{\nu +1},\theta}\left( { \omega   } \right)}  \over
		                         {   \gcontinuous\left( { \omega   } \right)}}     \right)^2  d\omega,
\end{eqnarray*}
the limit being a consequence of   $\filter(\delta \omega)  \rightarrow 1$   
(which can directly be seen from $\delta \gdelta(\delta \omega) -  \gcontinuous\left( { \omega  } \right) \rightarrow 0$  as $\delta  \rightarrow 0$, this difference being shown in fact $O(\delta^{2 \nu })$ in
  the Proof of Lemma~2.1 of Girard (2016)).

\bigskip

\noindent{\bf Proof of Lemma 2.2.} By Lemma A.1, since  
$\left\vert \left[\filter(\lambda)\right]^k  - \left[a_{\magn,\delta\theta}^*(\lambda)\right]^k  \right\vert <  k \left\vert \filter(\lambda)  - a_{\magn,\delta\theta}^*(\lambda)  \right\vert =O(\delta^{2 \nu} )$ 
one can replace, with an accuracy which will be sufficient since ${2 \nu} > {{2 \nu \over 2 \nu +1}}$, the filter by its un-aliased version in ${\int { \left[ \filter(\lambda) \right]^k   }  {\rm d}\lambda }$. On the other hand,  by the change of variable $s = \lambda /(\delta\theta)^{{2 \nu \over 2 \nu +1}}  $ and an application of the dominated convergence theorem,
one can find that 
$$    (\delta\theta)^{-{2 \nu \over 2 \nu +1}}    {\int_0^{2 \pi} { \left[a_{\magn,\delta\theta}^*(\lambda)\right]^k   }  {\rm d}\lambda } \rightarrow  \int _{0 }^{\infty }    {  (1+   (2 \pi \Cnu b)^{{-1}}  s^{2 \nu +1})^{-k}} {\rm d}s
=   (2 \pi \Cnu b)^{\frac{1}{{2 \nu +1}}} 
 \int _{0 }^{\infty }    {  (1+   s^{2 \nu +1})^{-k}} {\rm d}s . $$
The claimed equivalents are now obtained from a  known expression, for $k \in \{1,2  \}$, of  the latter classic integral.

\bigskip

\noindent{\bf Proof of Theorem 3.3.} 
It is convenient to apply   Theorem~1.58 of Kessler et al.~(2012) which is
a general asymptotic-consistency result for estimating equation.
In fact, our Theorem~3.1 exactly states that 
the conditions required to apply their Theorem~1.58 
are well fulfilled for, in their notations
 (except we set here $\overline \theta := (b,  \theta)^T$), 
 the two-component estimating equation
 $G_n\left( \overline \theta \right) := (  \inVarNoise \vecty^T  \vecty/n - b -1,    \inVarNoise {\rm CGEM}(\magn,\theta) - b )^T$.
Precisely, in their notations, the required conditions (i) and (ii) (resp.~(iii)) are immediate consequence of Part 1 (resp. Part 3) of Theorem~3.1. 

\bigskip

\noindent{\bf Proof of Part  1 of Theorem 4.1.}  We shall apply Theorem~1.60 
of Kessler et al.~(2012)
to the  two-component  equation $G_n\left( \overline \theta \right)$  defined above.
It is clear,  form the high differientiabily regularity, already mentioned, of  $\gdelta$ and its strict positivity, that $G_n$ is continuously differentiable over $B \times \Theta$.
Denoting by  ${ \partial  \over  \partial \overline \theta^T}  \overline G_n(b, \theta)$ the Jacobian matrix, 
 letting  $M_{\delta,b_0,\theta_0}( b, \theta)  :=$  
  $$ \left[
\begin{array}{cc}
  -1   &  0\\
{ {- 1  \over 2 \pi} 
  \left(  {
   \Bigint {   { \left\{   {   {[  \filter ]^2    -2 { \left(  [  \filter ]^2 - [  \filter ]^3  \right)    {   \left( {  {\magn_0 g_0}  \over  {\magn   \gdelta} } -1 \right)} } } } \right\}  }  }
}  \right)}
&    
{ {- \magn  \over 2 \pi} 
  \left(  {
   \Bigint { { \left\{   {   {[  \filter ]^2    + { \left(  2 [  \filter ]^3 - [  \filter ]^2  \right)    {   \left( {  {\magn_0 g_0}  \over  {\magn   \gdelta} } -1 \right)} } } } \right\}  }   \hdelta}
}  \right)}
\end{array}
\right]
, $$
 it is a consequence of Theorem~3.1 
that ${ \partial  \over  \partial \overline \theta^T}  \overline G_n(b, \theta)  \rightarrow M_{\delta,b_0,\theta_0}( b, \theta) $, uniformly over $B \times  \Theta$
and  $M_{\delta,b_0,\theta_0}( b_0, \theta_0) $ is invertible.
Furthermore, 
$
n^{1/2}
\overline G_n\left( \overline \theta_0 \right)
\xrightarrow[]{\cal D}
N\left( 
\left[
\begin{array}{ccc}
0 \\
0
\end{array}
\right],
 2 b_0^2
\left[
\begin{array}{cc}
 2  \pi \int   a_0^{-2}  g_0^2     &  1\\
 1 &       {1  \over 2 \pi}  \int   a_0^2
\end{array}
\right]
\right)
$
can be deduced
from classic time-series results (e.g.~in Azencott and Dacunha-Castelle (1986)) since the  required regularity conditions are well fulfilled  (as discussed in Section 3 and in Girard (2016))).
Then Theorem~1.60 
of Kessler et al.~(2012)  is applicable and it gives (after direct, albeit tedious, algebraic manipulations) that the asymptotic covariance matrix of $(\hat\magn_{\rm EV\vert \sigmaN},  \hat\theta_ {\rm GEV\vert \sigmaN})$ 
is $
{4  \pi
}
\left[
\begin{array}{cc}
 b_0^2   \int    a_0^{-2}  g_0^2     &  b_0    {   {  1 - \left(  \int   a_0^2 \int    a_0^{-2}  g_0^2   \right)    }\over   \int    a_0^2 h_0}\\
 b_0    {   {  1 - \left(  \int   a_0^2 \int   a_0^{-2}   g_0^2   \right)    }\over   \int   a_0^2  h_0}  &    \int   a_0^2   {   {   \left(  \int   a_0^2 \int  a_0^{-2}   g_0^2   \right)   -1  }\over   {\left\vert \int { a_0^2  h_0}\right\vert^{2} }}
\end{array}
\right].
$

\end{document}